\documentclass{article}

\newcommand{\Pre}{\textbf{PreO}}
\newcommand{\Pref}{\ensuremath{\textbf{PreO}_f}}
\newcommand{\Prelf}{\ensuremath{\textbf{PreO}_{lf}}}
\newcommand{\KFr}{\textbf{KFr}}
\newcommand{\CABA}{\textbf{CABA}}
\newcommand{\MA}{\textbf{MA}}
\newcommand{\Int}{\textbf{S4}}
\newcommand{\Intf}{\textbf{S4}_f}
\newcommand{\MAf}{\ensuremath{\textbf{MA}_f}}
\newcommand{\CAMA}{\ensuremath{\textbf{CAMA}_{\infty}}}
\newcommand{\SET}{\textbf{Set}}
\newcommand{\IndC}{\ensuremath{\rm{Ind}\-\mathcal{C}}}
\newcommand{\ProC}{\ensuremath{\rm{Pro}\-\mathcal{C}}}
\newcommand{\KFrf}{\ensuremath{\textbf{KFr}_{f}}}
\newcommand{\KFrlf}{\ensuremath{\textbf{KFr}_{lf}}}

\newcommand{\lra}{\longrightarrow}

\newcommand{\mC}{\mathcal{C}}

\RequirePackage{amsmath}

\RequirePackage[all]{xy}
\RequirePackage{hyperref}

\renewcommand{\iff}{\Leftrightarrow}

\newcommand{\dpll}[1]{{\sc DPLL}\xspace}
\newcommand{\Atm}{{\ensuremath{\mathit{Atoms}}}}

\RequirePackage{color}
\newcommand{\red}[1]{\textcolor{red}{#1}}

\newcommand{\COMMENT}[1]{}
\RequirePackage[normalem]{ulem}
\RequirePackage{ifthen}
\newcommand{\correction}[3][]{
  \textcolor{red}{\sout{#2}}\textcolor{blue}{#3}
  \ifthenelse{\equal{#1}{}}{}{\footnote{\red{#1}}}
}

\usepackage{arxiv}
\usepackage{amsthm, amsmath, amssymb, amsfonts , graphicx, multicol, array}
\usepackage[english]{babel}
\usepackage[utf8]{inputenc} 
\usepackage[T1]{fontenc}    
\usepackage{hyperref}       
\hypersetup{hypertexnames = false,
bookmarksdepth = 2,
bookmarksopen = true,
colorlinks,
linkcolor = black,
citecolor = black,
urlcolor = blue,
pdfstartview={XYZ null null 1},
}
\usepackage{url}            
\usepackage{booktabs}       
\usepackage{nicefrac}       
\usepackage{microtype}      
\usepackage{doi}
\usepackage{faktor}
\usepackage{tikz-cd}
\usepackage{comment}

\theoremstyle{plain}
\newtheorem{theorem}{Theorem}[section]
\newtheorem{lemma}[theorem]{Lemma}
\newtheorem{proposition}[theorem]{Proposition}
\newtheorem{corollary}[theorem]{Corollary}

\theoremstyle{definition}
\newtheorem{definition}[theorem]{Definition}
\newtheorem{example}[theorem]{Example}

\theoremstyle{remark}
\newtheorem{remark}[theorem]{Remark}

\title{Profiniteness, Monadicity and Universal Models in Modal Logic}

\author{{\hspace{1mm}Matteo De Berardinis}\thanks{Corresponding author} \\
	Department of Mathematics\\
	Università degli Studi di Milano\\
        Italy\\
	\And
	\href{https://orcid.org/0000-0001-6449-6883}{\includegraphics[scale=0.06]{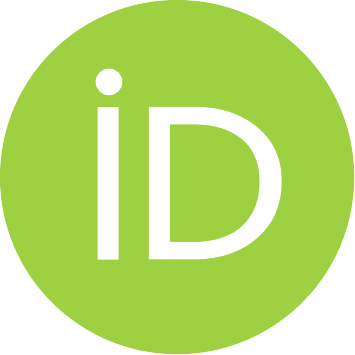}\hspace{1mm}Silvio Ghilardi} \\
	Department of Mathematics\\
	Università degli Studi di Milano\\
        Italy\\
}

\renewcommand{\shorttitle}{}

\hypersetup{
pdftitle={Profiniteness, Monadicity and Universal Models in Modal Logic},
pdfsubject={Modal Logic},
pdfauthor={Matteo De Berardinis, Silvio Ghilardi},
pdfkeywords={Monadic Functors, Profinite Algebras, Modal Algebras},
}

\begin{document}

\begin{abstract}
Taking inspiration from the monadicity of complete atomic Boolean algebras,
we prove that profinite modal algebras are monadic over $\bf Set$. While analyzing the monadic functor, we recover the universal model construction - a construction widely used in the modal logic literature for describing finitely generated free modal algebras and the essentially finite subframes of their canonical models.
\end{abstract}


\maketitle

\section{Introduction}\label{sec:intro}

Modal logic is an inspiring area for algebraic and coalgebraic techniques, as witnessed by a large literature from the last decades.
On one side, modal algebras (arising as Lindembaum algebras of  modal calculi) appear to be algebras for an endofunctor of the category of Boolean algebras; on the other side Kripke frames (the  main ingredient of the most popular semantics of modal logics) appear to be coalgebras for the covariant powerset endofunctor of the category of sets. In this paper, we shall investigate a different connection between modal logic and algebra/coalgebras constructions: this time, we shall take into consideration monads/co-monads and  algebras/co-algebras in the Eilenberg-Moore sense.

It is a basic result in elementary topos theory (see e.g.~\cite{Elephant}, Thm. A2.2.7) that the contravariant power set functor $P: \mathcal{E}^{op} \longrightarrow \mathcal{E}$ is monadic, thus making $\mathcal{E}^{op}$ equivalent to the category of the
algebras for the monad induced by the adjunction. In the case of $\mathcal{E}=\SET$, it is well-known that $\SET^{op}$ is equivalent to the category \CABA\ of complete atomic Boolean algebras, a category which then turns out to be monadic over $\SET$.

If one tries to reproduce the above result in the context of modal logic, one immediately
encounters seemingly unsolvable problems. If we let the category of Kripke frames and p-morphisms play the role of the category of  sets, we can easily identify (via Thomason duality~\cite{thomason}) the dual of this category as the category \CAMA\ of complete atomic completely-additive  modal algebras:
here a modal algebra is said to be `complete atomic completely-additive' if its underlying Boolean algebra is complete atomic and the possibility operator $\Diamond$  commutes with arbitrary Joins. However, the forgetful functor from \CAMA\ to the category of sets does not even have a left adjoint: in fact,  as we shall formally prove in Section~\ref{sec:Kfr} below, \CAMA\  does not  have an initial object and such an initial object should be carried over by a hypothetic left adjoint, because \SET\ has it.

However, this is only part of the story. In fact, if  we see complete atomic Boolean algebras as \emph{profinite} Boolean algebras (namely the Pro-Completion of the category of finite Boolean algebras) and then if we replace complete atomic completely additive modal algebras by profinite modal algebras,  the monadicity result can surprisingly be restored.

Monadicity over \SET\ is an important property: by definition, it means equivalence with the Eilenberg-Moore category of algebras for the monad induced by an adjunction. Monadicity can be seen as a generalized notion of `being algebraic': in fact, according to the well-known characterization~\cite{Bor,EV}, monadic categories over \SET\ share relevant properties (like exactness) with customary finitary varieties. The fact that
profinite modal algebras are monadic over \SET\ is a quite peculiar fact, because there is no reason why profinite algebras should be monadic (profinite sets are Stone spaces and the category of Stone spaces is not monadic over \SET\ because it is not exact).

Monadicity of profinite modal algebras has a \emph{coalgebraic} counterpart (also proved in the paper), namely comonadicity over \SET\ of the category of \emph{locally finite} Kripke frames.
The comonad induced by the adjunction between
locally finite Kripke frames and \SET\ is interesting by itself, because it turns out to be
strictly related to
the \emph{universal model} construction, namely a construction widely investigated in the modal logic literature for transitive modal systems~\cite{leo,valentin,be1,be2,ungh} in connection to
problems like
representation theorems,  atomicity of finitely generated free algebras, local finiteness of subvarieties, etc.

The paper is structured as follows: we first review basic facts concerning Ind- and Pro-completions, then we investigate relevant features of the category of locally finite Kripke frames leading to the proof of
the monadicity/comonadicity theorems. Finally, after showing that our monad is connected to the profinite completion of free algebras, we give an explicit description of it in the case of transitive modal systems.

As for prerequisites, the paper requires just basic category-theoretic background on monads, adjoints, limits and co-limits (see e.g. the textbook~\cite{CWM}).

\section{Review of Ind- and Pro-Completions}\label{sec:indpro}

In this Section we review some definitions and results from Chapter VI of~\cite{stone} about Ind- and Pro-Completions
 (we refer to that textbook for full proofs of the results stated here). The Ind-Completion of a category $\mC$ formally adds filtered colimits to $\mC$. One way of introducing it is as a  full subcategory of the presheaf category:


\begin{definition} For a small category $\mC$, we let  $\IndC$ be
 the full subcategory of $\mathbf{Set}^{\mathcal{C}^{\text{op}}}$ given by those functors that are filtered colimits of representable functors.
\end{definition}

We are only interested in the case where $\mC$ has finite colimits;
the following theorem summarizes the relevant properties of \IndC\ in  this hypothesis:

\begin{theorem}\label{thm:indchar}
Let $\mC$ be a small category with finite colimits. Then
\begin{itemize}
   \item[{\rm (i)}] \IndC\ is equivalent to $\text{Lex}(\mathcal{C}^{\text{op}},\SET)$, i.e. to the full subcategory of $\SET^{\mathcal{C}^{\text{op}}}$ given by the
 contravariant functors from $\mC$ to \SET\ turning finite colimits into finite limits;
    \item[{\rm (ii)}] \IndC\ has all filtered colimits and the embedding $\IndC \hookrightarrow \SET^{{\mC}^{\text{op}}}$ preserves them;
    \item[{\rm (iii)}] \IndC\ has  finite colimits and  the Yoneda embedding (restricted in the co\-do\-main)  $Y : \mC \hookrightarrow \IndC$ preserves them;
    \item[{\rm (iv)}] \IndC\  is complete, and the embedding $\IndC \hookrightarrow \SET^{{\mC}^{\text{op}}}$ preserves all small limits.
\end{itemize}
\end{theorem}

Notice that, as a consequence of (ii) and (iii), we have that \IndC\ is co-complete;
actually it is a co-completion of $\mC$ in the sense of the following:

\begin{definition}
Given a small category $\mC$, a \textit{cocompletion} of $\mC$ is a full embedding $F : \mC \rightarrow \mathcal{D}$ into a cocomplete category $\mathcal{D}$ s.t. every object of $\mathcal{D}$ is a colimit of objects in the image of $F$.
\end{definition}

 We recall the notion of finite presentability.

\begin{definition}
An object $X$ of a (locally small) category with filtered colimits $\mathcal{D}$ is said to be \textit{finitely-presentable} (in $\mathcal{D}$) if the functor $\text{Hom}_{\mathcal{D}}(X,-) : \mathcal{D} \rightarrow \SET$ preserves filtered colimits.
\end{definition}

The following theorem (to be used in the sequel) characterizes \IndC\ up to equivalence as that cocompletion of $\mC$ for which the embedding functor preserves finite colimits and sends the objects of $\mC$ to finitely-presentable objects.

\begin{theorem}\label{teo:charindcompl}
Let $\mC$ be a small category with finite colimits, and $Z : \mC \hookrightarrow \mathcal{D}$ a full embedding of $\mC$ in a cocomplete category $\mathcal{D}$. Then
\begin{itemize}
    \item[(i)] if the objects in the image of $Z$ are finitely-presentable in $\mathcal{D}$, $Z$ extends to a full embedding $\hat Z : \IndC \hookrightarrow \mathcal{D}$;
    \item[(ii)] if in addition $Z : \mC \hookrightarrow \mathcal{D}$ is a cocompletion of $\mC$ and $Z$ preserves finite colimits, $\hat Z$ is an equivalence.
\end{itemize}
\end{theorem}

The dual of the notion of Ind-Completion is the notion of profinite completion:

\begin{definition}\label{def:pro}
If $\mathcal{C}$ is a small category with finite limits, we let $\ProC$ to be the category  defined as $(\rm{Ind}\-(\mathcal{C}^{op}))^{op}$ (that is, $\ProC$ is dual to the category of covariant finite-limit preserving
functors from $\mC$ to \SET).
\end{definition}

We conclude the section with a couple of examples (also taken from~\cite{stone}), which are relevant for  our paper.

\begin{example}
Given a variety (i.e. an equationally defined class of algebras) $\mathcal{V}$, one can apply Theorem~\ref{teo:charindcompl} to show that $\mathcal{V}$ is the Ind-Completion of the category of finitely presented $\mathcal{V}$-algebras. In particular, Boolean algebras are the Ind-Completion of the category of finite Boolean algebras.\footnote{
In general, the notion of a finite algebra and of a finitely presented algebra do not coincide; they coincide for locally finite varieties, where finitely generated free algebras are finite (this is the case of Boolean algebras, but not of modal algebras).
} Finite sets are dual to finite Boolean algebras and the category of Stone spaces and continuous maps is dual to the category of Boolean algebras; thus one can conclude that the category of Stone spaces is the Pro-Completion of the category of finite sets.
\end{example}

\begin{example}
The category of sets is the variety defined by the empty set of equations over the empty signature; thus \SET\ is the Ind-Completion of the category $\SET_{fin}$ of finite sets. Since  $\SET_{fin}$ is dual to the category of finite Boolean algebras and the category \CABA\ of complete atomic Boolean algebras is dual to \SET\ (see Theorem~\ref{thm:Tarski} below), we conclude that \CABA\ is the Pro-Completion of the category of finite Boolean algebras (otherwise said, the complete atomic Boolean algebras are precisely  the profinite Boolean algebras).
\end{example}


\section{Modal Algebras and Kripke Frames}\label{sec:Kfr}

In this section we supply some  algebraic background and recall Tarski and Thomason dualities (proofs are folklore, we give them in the Appendix); then we investigate colimit constructions in the category of Kripke frames.

Modal algebras are standard algebraic semantics for the minimum normal modal logic $K$; varieties of modal algebras are in one-to-one correspondence with propositional normal modal logics. The reader is referred to the textboox~\cite{misha} for comprehensive information on modal logic, we just recall here the definitions that are relevant for the paper.

A \emph{modal algebra} $(B, \Diamond)$ is a Boolean algebra $B$ endowed with a hemimorphism, i.e.
with a finite-join preserving operator $\Diamond: B \rightarrow B$:
$$
\Diamond (a\vee b)= \Diamond a \vee \Diamond b, \qquad \Diamond \bot =\bot~.
$$
Below we shall use the letter $B$ to mean both a modal algebra and its support set; the $\Diamond$ operator is called the `possibility' operator of the modal algebra and its dual (defined as $\Box a:= \neg \Diamond \neg a$) is called the necessity operator of the modal algebra.
Modal algebras and Boolean morphisms preserving $\Diamond$ form the category $\MA$;
we let $\MAf$ be the full subcategory of $\MA$ formed by the finite modal algebras.

On the semantic side, we have the notions of a \emph{Kripke frame} and of a p-morphism between Kripke frames.

\begin{definition}
    A \emph{Kripke frame} is a directed graph, i.e. a pair $(P,R)$ given by a set $P$ and a binary relation $R$ on it.
\end{definition}

\begin{definition}
    Given two Kripke frames $(P,R)$ and $(Q,S)$, a \emph{p-morphism} from $(P,R)$ to $(Q,S)$ is a function $f : P \rightarrow Q$ with the following properties:
    \begin{itemize}
        \item[$\bullet$]
        (\textit{stability})~
        {for $p,p' \in P$, $p R p' \implies f(p) S f(p')$;}
        \item[$\bullet$]
        (\textit{openness})
        {for $q' \in Q$, $p \in P$, $f(p) S q' \implies \exists p' \in P \text{ s.t. } p R p' ~\&~
        f(p') = q'$.}
    \end{itemize}
\end{definition}

We define \KFr\ as the category having Kripke frames as objects and p-morphisms as morphisms; \textbf{DGrph} will be the category having the same objects of \KFr\ and stable functions as morphisms.\footnote{The name `stable function' comes from recent modal logic literature, see e.g.~\cite{BBI}.} Sometimes, when referring to an object of \KFr\ (or \textbf{DGrph}), we will omit to write the binary relation.

We now recall Tarski and Thomason dualities, relating algebraic and semantic notions. The following proposition is well-known:

\begin{proposition}\label{prop:caba}
 The following conditions are equivalent for a Boolean algebra $B$:
  \begin{itemize}
   \item[{\rm (i)}] $B$ is complete and atomic;
   \item[{\rm (ii)}] $B$ is isomorphic to a powerset Boolean algebra;
   \item[{\rm (iii)}] $B$ is complete and satisfies the infinitary distributive law:
     \begin{align*}
      \bigwedge_{i \in I} \bigvee X_i = \bigvee_{f \in \prod_{i \in I} X_i} \bigwedge_{i \in I} f(i)
     \end{align*}
     for every family $\{X_i\}_{i \in I}$ of subsets of $B$.
  \end{itemize}
\end{proposition}

Complete atomic Boolean algebras and Boolean morphisms preserving arbitrary Joins and Meets form the category $\CABA$. A modal algebra $(B, \Diamond)$ is said to be \emph{completely additive complete atomic} iff $B$ is complete atomic as a Boolean algebra and $\Diamond$ preserves arbitrary Joins. Completely additive complete atomic modal algebras, endowed with Boolean morphisms preserving $\Diamond$ as well as arbitrary Joins and Meets form the category $\CAMA$.

\begin{theorem}\label{thm:Tarski}
[Tarski duality] \CABA\ is dual to \SET.
\end{theorem}

\begin{proof}
 See the Appendix. We only recall how the duality works. From \SET\ to \CABA, the duality functor maps a set to its powerset and a function to the inverse image Boolean morphism. From \CABA\ to \SET, the duality functor associates with a complete atomic Boolean algebra the set of its atoms and to a Boolean morphism preserving arbitrary Joins and Meets its posetal left adjoint restricted to atoms (in the domain and in the codomain).
\end{proof}

\begin{theorem}\label{thm:Thomason}
 [Thomason duality] \CAMA\ is dual to \KFr.
\end{theorem}

\begin{proof}
 The duality functors extend Tarski duality functors. On the Boolean algebra side, the duality functors make the set of atoms of a Boolean algebra into a Kripke frame by setting $aRb$ iff $a\leq \Diamond b$.
 On the side of a Kripke frame $(P, R)$, the Boolean algebra $\mathcal{P}(P)$ is turned into a modal algebra by setting $\Diamond X:= \{p\in P \mid \exists p'\; (pRp' \;\& \, p'\in X)\}$.
\end{proof}

We now recall standard semantic constructions on Kripke frames.

\begin{definition}
A subset $G \subseteq P$ is called \emph{generated subframe} (of the Kripke frame $(P,R)$) if, for all $p, p' \in P$,
\begin{align*}
p \in G ~\&~ p R p' \implies p' \in G~~~.
\end{align*}
\end{definition}

\begin{definition}
    If there is a p-morphism $(P,R) \rightarrow (Q,S)$ that is surjective as a function, $(Q,S)$ is said to be a \emph{p-morphic image} (of $(P,R)$).
\end{definition}

The above constructions can be  exploited categorically:

\begin{proposition}\label{prop:cocomp}
\textbf{KFr} is cocomplete. Moreover, the forgetful functors $\textbf{KFr} \rightarrow \textbf{Set}$ and $\textbf{KFr} \rightarrow \textbf{DGrph}$ preserve colimits.
\end{proposition}
\begin{proof}
We compute coproducts and coequalizers.

Coproducts are easy:
    given a family of objects $\{P_i\}_{i \in I}$ in \KFr, the coproduct is the disjoint union $\coprod P_i$ of the sets $P_i$ endowed with the relation $\coprod R_i$ given by the union of the $R_i$'s; by definition, the inclusions $\iota_j : P_j \rightarrow \coprod P_i$ define p-morphisms $(P_j,R_j) \rightarrow (\coprod P_i, \coprod R_i)$ for all $j \in I$. The universal property of coproducts is easily checked.

    Coequalizers are more interesting.
    Given a pair of parallel morphisms $f, g: P \rightrightarrows Q$ in \KFr, let's consider the quotient of $Q$ over $\sim$, the equivalence relation generated by $\{(f(p),g(p)) \, | \, p \in P\}$ (i.e. the \emph{smallest} equivalence relation containing it); we will use $[x]$ to indicate the class of $x$ wrt the equivalence relation $\sim$. On $\faktor{Q}{\sim}$, we consider the following binary relation: for $x,y \in Q$,
    \begin{align*}
    [x]\overline{S}[y] \iff \exists \tilde{y} \in Q \text{ s.t. } \tilde{y} \sim y \text{ and } x S \tilde{y}
    \end{align*}

    The relation $\overline{S}$ is well defined (we have to prove the independence for the choice of the element in the class $[x]$). For this purpose, let's define another equivalence relation on $Q$: for $x,x' \in Q$,
    \begin{align*}
    x \equiv x' \iff \bigcup_{\tilde{y} \text{ s.t. } xS\tilde{y}} [\tilde{y}] = \bigcup_{\tilde{y}' \text{ s.t. } x'S\tilde{y}'} [\tilde{y}']
    \end{align*}
    Fixing $p \in P$ and $\tilde{y} \in Q$ s.t. $f(p)S\tilde{y}$, there exists $\tilde{p} \in P$ s.t. $pR\tilde{p}$ and $f(\tilde{p})=\tilde{y}$; if we pick some $y \in Q$ s.t. $y \sim \tilde{y}$, then $y \sim f(\tilde{p}) \sim g(\tilde{p})$ and $g(p)Sg(\tilde{p})$. This proves that
    \begin{align*}
    \bigcup_{\tilde{y} \text{ s.t. } f(p)S\tilde{y}} [\tilde{y}] = \bigcup_{\tilde{y}' \text{ s.t. } g(p)S\tilde{y}'} [\tilde{y}']
    \end{align*}
    for all $p \in P$ (we proved $\subseteq$, but the other inclusion is analogous), and this means that $\{(f(p),g(p)) \, | \, p \in P\} \subseteq\, \equiv$. Now, by definition of $\sim$, we can conclude that $\sim \,\subseteq\, \equiv$, i.e. that $\overline{S}$ is well defined.

    By definition, the projection $\pi : Q \rightarrow \faktor{Q}{\sim}$ defines a p-morphism $(Q,S) \rightarrow \left(\faktor{Q}{\sim},\overline{S}\right)$.

    Now, given a set $A$, together with a function $h : Q \rightarrow A$ s.t. $h\circ f = h\circ g$, there exists a unique function $\phi : \faktor{Q}{\sim} \rightarrow A$ s.t. $\phi\circ \pi = h$ (it must send $[x]$ into $h(x)$).

    If $A$ is equipped with a binary relation $T$ and $h$ is stable, then $\phi$ is stable: if $[x]\overline{S}[y]$ in $\faktor{Q}{\sim}$, let $\tilde{y} \in Q$ be s.t. $\tilde{y} \sim y$ and $xS\tilde{y}$; then $\phi([x]) = h(x) T h(\tilde{y}) = \phi([\tilde{y}]) = \phi([y])$.

    If $h$ is open, then $\phi$ is open: if $\phi([x]) T a$ in $A$, then, since $\phi([x])=h(x)$, there exists $y \in Q$ s.t. $xSy$ and $h(y) = a$, hence $[x]\overline{S}[y]$ and $\phi([y])=a$.
\end{proof}
These computations show not only that  \KFr\ has all colimits, but also that the forgetful functors
$\KFr \lra \SET$ and  $\KFr \lra \textbf{DGrph}$
preserve those colimits.
However, these functors do not have right adjoints; in fact
\KFr\ (and consequently also \CAMA) is not a good category to consider, as it lacks  minimal properties:

\begin{proposition}
 \KFr\ does not have a terminal object.
\end{proposition}

\begin{proof}
 Consider an ordinal $(\alpha, >)$ seen as a Kripke frame. Notice that any p-morphism $f: (\alpha, >) \lra (P, R)$ must be injective: this is becase if we have $f(p_0)=f(p'_0):=q$ for $p_0>p'_0$, then we get $qRq$ by stability.
 But then, by openness, from $f(p_0)R q$ we see that there is $p_1$ with $p_0 > p_1$ such that $f(p_1)=q$.
 Continuing in this way, we produce an infinite descending chain $p_0>p_1>p_2\cdots$, which cannot
 exist.
 Thus a hypothetical terminal object in \KFr\ must contain a copy of every ordinal: this  is too large to 
 be a set.
\end{proof}

Equalizers however are easily computed.
Let's consider a pair of parallel morphisms $f, g : P \rightrightarrows Q$ in \textbf{KFr}; we set
\begin{align*}
E_{fg} := \{p \in P \; \mid \; \forall p'\; (pRp' \Rightarrow f(p') = g(p')) \; \}
\end{align*}
$E_{fg} \subseteq P$ is obviously a generated subframe and it is the biggest generated subframe contained in $\{p \in P \; | \; f(p)=g(p)\} \subseteq P$.

\begin{lemma}\label{lem:eq}
 The inclusion $\iota : E_{fg} \rightarrow P$ is the equalizer of $f$ and $g$ in \KFr.
\end{lemma}

\begin{proof}
 $\iota$ is  clearly a morphism in \KFr. To check that $\iota$ is the equalizer of $f$ and $g$, consider  another morphism $h : A \rightarrow P$ (let's say $A$ equipped with the binary relation $T$) s.t. $f\circ h = g\circ h$. If we fix $x \in A$ and $p'\in P$ such that  $h(x)R p'$, we can find $x' \in A$ s.t.
 $xTx'$ and $h(x')=p'$, implying $f(p') = f(h(x')) = g(h(x')) = g(p')$; this proves that $\text{Im}h \subseteq E_{fg}$, hence we have a unique factorization in \KFr\ of $h$ trough $\iota$.
\end{proof}

\begin{remark}\label{rem:eq} Hence
the forgetful functors $\KFr \rightarrow \SET$ and $\KFr \rightarrow \textbf{DGrph}$ do not preserve equalizers. However, they preserve equalizers of those pairs $f, g : P \rightrightarrows Q$ for which the set $\{p \in P \; | \; f(p)=g(p)\}$ is a generated subframe of $P$.
\end{remark}

\section{Locally finite Kripke frames}\label{sec:kfrlf}

There exists a nice full subcategory of \KFr\ worth investigating.
Let $(P, R)$ be a Kripke frame; we indicate with $R^*$ the reflexive-transitive closure of $R$, namely
for $p,q\in P$ we define
\begin{align*}
p R^* p'~~~~~ \iff~~~~~ \exists n \in \mathbb{N}\; \exists p_0,\dots,p_n \in P \text{ s.t. } p=p_0 R \dots R p_n=p'
\end{align*}
For $p \in P$, we write $R(p)$ and $R^*(p)$ for $\{q\in P\mid pRq\}$ and $\{q\in P\mid pR^*q\}$, respectively.
Notice that $R^*(p)$ is the smallest generated subframe of $P$ containing $p$.

\begin{definition}
    A Kripke frame $(P,R)$ is said to be \emph{locally-finite} if $R^*(p)$ is finite for all $p \in P$.
\end{definition}

We let \KFrlf\ be the full subcategory of \KFr\ given by the locally finite Kripke frames and we let
\KFrf\ be the full subcategory of \KFr\ given by the  finite Kripke frames.

 \KFrlf\ and \KFrf\ are closed under generated subframes, p-morphic images and finite disjoint unions  (\KFrlf\ also under infinite disjoint unions); as a consequence, from the proof of Proposition~\ref{prop:cocomp}, we obtain:

\begin{proposition}\label{prop:cocomplete1}
\KFrlf\ is cocomplete and \KFrf\ has finite colimits. Moreover, the embedding $\KFrlf\ \hookrightarrow \KFr$ preserves colimits and the embedding $\KFrf\ \hookrightarrow \KFrlf$ preserves finite colimits.
\end{proposition}

We can now apply the chracterization of Ind-Completions given by Theorem~\ref{teo:charindcompl}:

\begin{theorem}\label{thm:indcfin}
     \KFrlf\
     is the Ind-Completion of \KFrf. Moreover, \KFrlf\ coincides with the full subcategory of \KFr\ given by those objects that are expressible as filtered colimits in \KFr\ of diagrams in \KFrf.
\end{theorem}
\begin{proof}
Notice that
the category \KFrf\ is not small, but it is essentially small, i.e. it is equivalent to a small one.
We apply Theorem~\ref{teo:charindcompl}.

$\KFrf \hookrightarrow \KFrlf$ is a full embedding and in Proposition~\ref{prop:cocomplete1}, we saw that \KFrlf\ is cocomplete.

We show that every object of \KFrlf\  is a colimit of objects in \KFrf. Given $(P,R)$ in \KFrlf, let's consider the diagram in $\KFrf$ given by the finite Kripke frames $R^*(p)$, with $p \in P$, and the inclusions $R^*(p) \hookrightarrow R^*(p')$ (whenever $pR^*p'$); it's straightforward to see that $P$ (with the inclusions $R^*(p) \rightarrow P$) is the colimit of this diagram in \KFrlf\ (we might as well have used the diagram given by the generated subframes of $P$ that have a finite underlying set - in this case, the colimit would have been  filtered). This argument proves also the second assertion of the theorem, since $ \KFrlf\hookrightarrow \KFr$ preserves colimits.

By Proposition~\ref{prop:cocomplete1}, we know that \KFrf\  has all finite colimits and that
$\KFrf\ \hookrightarrow \KFrlf$ preserves them.
It remains to prove that the objects of \KFrf\ are finitely-presentable in \KFrlf.

Let's consider a $Q$ in \KFrlf, colimit of some filtered diagram $I \rightarrow \KFrlf$, with morphisms $\varphi_i : Q_i \rightarrow Q$.
%
    Since the forgetful functor $\KFrlf \rightarrow \SET$ (composite of the embedding $\KFrlf \hookrightarrow \KFr$ and of the forgetful $\KFr \rightarrow \SET$) preserves colimits, $Q$ is isomorphic, as a set, to a quotient of $\coprod_i Q_i$ and, for $i,j \in I$, $x \in Q_i$ and $x' \in Q_j$, we have that $\varphi_i(x) = \varphi_{j}(x')$ iff there exist $i \rightarrow k$ and $j \rightarrow k$ in $I$ s.t. the induced $Q_i \rightarrow Q_k$ and $Q_j \rightarrow Q_k$ send $x$ and $x'$ into the same element.

We have to prove that the induced morphism
$$\underrightarrow{\textit{lim}}_i \text{Hom}_{\KFrlf}(P,Q_i) \rightarrow \text{Hom}_{\KFrlf}(P,\underrightarrow{\textit{lim}}_i Q_i)$$
is a bijection for all objects $P$ in $\KFrf$.

Let's consider $f : P \rightarrow Q \cong \underrightarrow{\textit{lim}}_i Q_i$ in $\KFrlf$; if we see this  as a morphism in $\SET$, since the forgetful functor $\KFrlf \rightarrow \SET$ preserves colimits and finite sets are finitely-presentable in $\SET$, there exists a unique factorization (up to equivalence in $\underrightarrow{\textit{lim}}_i \text{Hom}_{\SET}(P,Q_i)$)
\begin{center}
\begin{tikzcd}
                                          & Q_i \arrow[d, "\varphi_i"] \\
P \arrow[r, "f"'] \arrow[ru, "\tilde{f}"] & Q
\end{tikzcd}
\end{center}
In general, $\tilde{f}$ is just a function, but it is stable for a suitable choice of the index $i$.
In fact,\footnote{
Stability of $\tilde f$ is alternatively ensured by the fact that finite graphs are finitely presentable in $\textbf{DGrph}$ and by the fact that the inclusion $\KFrlf \hookrightarrow \textbf{DGrph}$ preserves colimits (see Propositions~\ref{prop:cocomp} and~\ref{prop:cocomplete1}).
}
if we have $p,p' \in P$ s.t. $p R p'$, then ($f$ is a p-morphism) $\varphi_i(\tilde{f}(p)) = f(p) S f(p') = \varphi_i(\tilde{f}(p'))$; since $\varphi_i$ is a p-morphism, there exists $x' \in Q_i$ s.t. $\tilde{f}(p) S_i x'$ and $\varphi_i(x') = \varphi_i(\tilde{f}(p'))$. Now, using the previous remark and the fact that any two parallel arrows $\alpha, \beta : i \rightrightarrows j$ in $I$ can be coequalized by some $\gamma : j \rightarrow k$ in $I$ (by the properties of a filtered category), we have a morphism $i \rightarrow k$ in $I$ s.t. the induced $Q_i \rightarrow Q_k$ sends $x'$ and $\tilde{f}(p')$ into the same element; this means that, up to composition with $Q_i \rightarrow Q_k$, we have $\tilde{f}(p) S_i \tilde{f}(p')$. Moreover, if $\tilde{f}$ already has the stability property for a certain pair $pRp'$ in $P$, the composition of $\tilde{f}$ with $Q_i \rightarrow Q_k$ (which is a p-morphism) does too. Being $P$ finite as a set, we can assume that $\tilde{f}$ is a stable function (there are a finite number of pairs $pRp'$ in $P$ and, after each composition with $Q_i \rightarrow Q_k$, the number of those pairs for which $\tilde{f}$ doesn't have the stability property is strictly lower).

Let now prove openness of $\tilde{f}$;
if we have $\tilde{f}(p) S_i x'$ in $Q_i$, then ($\varphi_i$ is a p-morphism) $f(p) = \varphi_i(\tilde{f}(p)) S \varphi_i(x')$, hence ($f$ is a p-morphism) there exists $p' \in P$ s.t. $pRp'$ and $\varphi_i(\tilde{f}(p'))=f(p')=\varphi_i(x')$;
if we compose with  $Q_i \buildrel Q_d\over\rightarrow Q_k$ (induced by a suitable $d:i \rightarrow k$ in $I$), we have that $Q_d(\tilde{f}(p')) = Q_d(x')$. We proved that

(*) ``for all $p\in P$, for all $x'$ such that  $\tilde{f}(p) S_i x'$, there are  $p' \in P$ and
$d:i \rightarrow k$ in $I$ s.t. $pRp'$ and $Q_d(\tilde{f}(p')) = Q_d(x')$.''

\noindent
Since $P$ is finite and $Q_i$ is locally finite, there are finitely many such pairs $(p, x')$ and consequently, being $I$ filtered, we can take the same $d$ for all such pairs. The composition
$Q_d\circ \tilde{f}$ now fits our purposes, because if we have $Q_d(\tilde{f}(p')) S_k y$ for some $y\in Q_k$, then
(as $Q_d$ is a p-morphism) there is $x'$ such that $\tilde{f}(p) S_i x'$ and $Q_d(x')=y$: applying (*) to the pair $(p,x')$ we get $p' \in P$ such that $pRp'$ and $Q_d(\tilde{f}(p')) = Q_d(x')=y$, proving that
$Q_d\circ \tilde{f}$ is open (otherwise said, $\tilde f:P \lra Q_i$ itself is open, for a suitable choice of the index $i$).
\end{proof}

If we now recall Theorem~\ref{thm:indchar}, we immediately obtain:
\begin{corollary}
 \KFrlf\ is equivalent to $\text{Lex}(\KFrf^{\text{op}},\SET)$; moreover, \KFrlf\ is complete and the embedding $\KFrlf\hookrightarrow \SET^{{\KFrf}^{\text{op}}}$ preserves all limits.
\end{corollary}

It should be noticed that, however, limits are hard to be computed directly in \KFrlf. The problem are products (for equalizers, it is easy to see that Lemma~\ref{lem:eq} holds for \KFrlf\ too).
The  idea to use products in the underlying category $\textbf{DGrph}$ is wrong (projections from the  $\textbf{DGrph}$-product are p-morphisms but the obvious candidates for universal maps into the $\textbf{DGrph}$-product need not be open).  Limits are easily computed in $\text{Lex}(\KFrf^{\text{op}},\SET)$, but the equivalence $\KFrlf \simeq \text{Lex}(\KFrf^{\text{op}},\SET)$ is non trivial.

On the contrary, colimits are easy to compute directly in \KFrlf, but are rather involved if computed in $\text{Lex}(\KFrf^{\text{op}},\SET)$.

A dual characterization for profinite modal algebras can be obtained immediately from Definition~\ref{def:pro} and from the fact that $\MAf\simeq (\KFrf)^{op}$:

\begin{theorem}\label{thm:proMAf}
 Let \MAf\ be the full subcategory of \MA\ given by the finite modal algebras. Then
 ${\rm Pro}\-\MAf$ is dual to \KFrlf.
\end{theorem}

An analogous result (but limited to  objects) characterizing profinite Heyting algebras  was obtained in~\cite{GuramNick}.

\section{Monadicity of Profinite Modal Algebras}

Now, let's go back to the forgetful functors $\KFrlf \rightarrow \SET$ and $\KFrlf \rightarrow \textbf{DGrph}$ (we will call both of them $U$); we already observed that they preserve all (small) colimits (being the composites of
the inclusion
$\KFrlf \hookrightarrow \KFr$ with the corresponding forgetful functors), so we wonder if they have right adjoints.

\begin{theorem}\label{thm:adjoint}
 The forgetful functors $\KFrlf \rightarrow \SET$ and $\KFrlf \rightarrow \textbf{DGrph}$ have a right adjoint.
\end{theorem}

\begin{proof}
We use the classical `Special Adjoint Functor Theorem' (SAFT) result~\cite{CWM}, in the dual version because we are looking for a right adjoint.
We check that the conditions of the dual SAFT are satisfied by $U$:
\begin{itemize}
    \item[$\bullet$] \KFrlf, \SET\ and $\textbf{DGrph}$ have small hom-sets;
    \item[$\bullet$] \KFrlf\ has all small colimits and $U$ preserves them (as proved in Propositions~\ref{prop:cocomplete1} and
    \ref{prop:cocomp});
    \item[$\bullet$] \KFrlf\ is well-copowered: this is clear because,  as shown in Lemma~\ref{lem:epi} below, the epimorphisms of \KFrlf\ coincide with the p-morphisms  that are surjective as functions.
    \item[$\bullet$] \KFrlf\ has a small generating set: in fact the objects coming from the essentially small subcategory  \KFrf\ are a generating set for \KFrlf, because  all objects of \KFrlf\  are colimits of some diagram from \KFrf.
\end{itemize}
\end{proof}

\begin{lemma}\label{lem:epi}
In \KFr, \KFrf\ and \KFrlf\,  epimorphisms coincide with the  morphisms that are surjective as functions.
\end{lemma}

\begin{proof}
 Let   $f : (P, R) \rightarrow (Q, S)$ be a morphism in one of the above categories and let's consider the disjoint union $A := Q + (Q \setminus \text{Im}f)$, with (injective) functions $\iota_1 : Q \rightarrow A$ (which sends every element into its copy in the first component of $A$) and $\iota_2 : Q \rightarrow A$ (which sends the elements of $\text{Im}f$ into their copy in the first component of $A$ and the elements of $Q \setminus \text{Im}f$ into their copy in the second component of $A$). On $A$, let's consider the following binary relation: for $a, b \in A$
\begin{align*}
a T b \iff \begin{cases}
                {\makebox[10cm][l]{$a=\iota_1(x), b=\iota_1(y), x,y \in Q \text{ and } xSy$} \text{or}}\\
                {\makebox[10cm][l]{$a=\iota_2(x), b=\iota_2(y), x,y \in Q \setminus \text{Im}f \text{ and } xSy$} \text{or}}\\
                {\makebox[10cm][l]{$a=\iota_2(x), b=\iota_1(y), x \in Q \setminus \text{Im}f, y \in \text{Im}f \text{ and } xSy$}}
             \end{cases}
\end{align*}
With this relation, $\iota_1$ and $\iota_2$ define p-morphisms (one could check that $\iota_1, \iota_2 : Q \rightrightarrows A$ is in fact the cokernel pair of $f : P \rightarrow Q$). Clearly $\iota_1\circ f=\iota_2\circ f$;
if $f$ is epi, then $\iota_1=\iota_2$ which can be only in case $Im(f)=Q$. Vice versa, it is obvious that surjective p-morphisms are epi.
\end{proof}

Putting together Theorem~\ref{thm:proMAf} and Tarski duality, we obtain:

\begin{theorem}
 The forgetful functor ${\rm Pro}\-\MAf \lra \CABA$ has a left adjoint.
\end{theorem}

\begin{proof}
 Just recall that Tarski duality associates to a set $X$ the complete atomic Boolean algebra $\mathcal{P}(X)$, so that the composite functor  ${\rm Pro}\-\MAf \simeq \KFrlf^{op}\buildrel{U^{op}}\over \lra \SET^{op} \simeq \CABA$ is the functor associating with a profinite modal algebra its underlying
 Boolean algebra (up to a natural bijection).
\end{proof}

\begin{theorem}\label{thm:comon}
The forgetful functor $U:\KFrlf \rightarrow \SET$ is co-monadic.
\end{theorem}

\begin{proof}
 We check the hypotheses of  the CTT (`Crude Tripleability Theorem', see~\cite{barr}, Section 3.5 or~\cite{Elephant}, Theorem A1.1.2), in the dual context because we want to prove comonadicity.
 The hypotheses are the following ones:
\begin{itemize}
    \item[$\bullet$] $U$ has right adjoint (see the previous theorem).
    \item[$\bullet$] $U$ reflects isomorphisms: in fact, it is easy to see that a morphism $f : (P,R) \rightarrow (Q,S)$ in \KFr\ (and in \KFrlf, \KFrf\ too) is an isomorphism  iff it is a bijection as a function.\footnote{Notice that this is not true for morphisms in \textbf{DGrph} and in fact the forgetful functor $\textbf{DGrph} \lra \SET$ is not comonadic.
    }
    \item[$\bullet$] \KFrlf\ has and $U$ preserves  equalizers of coreflexive pairs. In fact,  in $\KFrlf$ equalizers exist and are computed as in \KFr\ via Lemma~\ref{lem:eq}.  If we consider a
    a diagram
    \begin{center}
    \begin{tikzcd}
    (P,R) \arrow[r, "f", shift left=1.5] \arrow[r, "g"', shift right=1.5] & (Q,S) \arrow[l, "t"', bend right, shift right=3]
    \end{tikzcd}
    \end{center}
    s.t. $t\circ f=1_P$ and $t\circ g=1_P$, the set $\{p \in P \; | \; f(p)=g(p)\}$ is a generated subframe of $P$: if $p \in P$ is s.t. $f(p)=g(p)$ and $p R p'$, then $g(p) = f(p) S f(p')$, hence there exists $p_1 \in P$ s.t. $p R p_1$ and $g(p_1)=f(p')$; necessarily, $p_1=tg(p_1)=tf(p')=p'$, so we can conclude that $f(p') = g(p')$. Thus, by Remark~\ref{rem:eq},
     we conclude that $U$ preserves equalizers of coreflexive pairs.
\end{itemize}
\end{proof}

In general, monadic functors do not compose, however the composition of two monadic functors is monadic in case the first component of the composition satisfies CTT, see~\cite{barr}, Section 3.5. Thus, from the well-known fact that
the forgetful functor $\CABA \lra \SET$ is monadic
(\cite{Elephant}, Thm. A2.2.7)
and from Theorems~\ref{thm:proMAf} and~\ref{thm:comon}, we conclude our main result:

\begin{theorem}
 The forgetful functors ${\rm Pro}\-\MAf \lra \CABA$ and ${\rm Pro}\-\MAf \lra \SET$ are both monadic.
\end{theorem}

We now want to better investigate the left adjoint to the forgetful functor ${\rm Pro}\-\MAf \lra \SET$ and relate it to the profinite completions of free modal algebras.

\begin{lemma}
 The faithful inclusion $\CAMA\hookrightarrow \MA$ prerserves limits.
\end{lemma}
 \begin{proof}
  This is due to the fact that limits in \CAMA\ are computed as in \SET, because so happens in the underlying category of complete atomic Boolean algebras (see the infinitary equational axiomatization given in Proposition~\ref{prop:caba}(iii)).
 \end{proof}

 \begin{theorem}
  The inclusion functor $\iota:{\rm Pro}\MAf \hookrightarrow \MA$ has a left adjoint ${\rm PrC}: \MA \lra {\rm Pro}\MAf$.
 \end{theorem}

 \begin{proof}
Notice that $\iota$ is the composition of two inclusions
$${\rm Pro}\MAf\hookrightarrow \CAMA \hookrightarrow \MA,
$$
both of which are limit-preserving: the second inclusion preserves limits by the previous lemma and the first inclusion preserves limits because it is dual to the inclusion $\KFrlf \hookrightarrow \KFr$ which is known to preseve colimits by Proposition~\ref{prop:cocomplete1}.  Moreover \MA\ has small hom-sets and the remaining conditions for SAFT are satisfied by {\rm Pro}\MAf\ (we checked the duals of such conditions for \KFrlf\ in the proof of Theorem~\ref{thm:adjoint}).
 \end{proof}

The functor ${\rm PrC}$ is called the `profinite completion functor'; the reason for this terminology will be clear from the characterization of the unity of the above adjointness given below.

Let us indicate with $\iota_f$ the restrition of the functor $\iota:{\rm Pro}\MAf \hookrightarrow \MA$
to the subcategory \MAf\ in its domain and
let us consider, for each object $A$ in \MA, the essentially small comma category $A \downarrow \iota_f$  and the forgetful functor $\pi : A \downarrow \iota_f \rightarrow \MAf  \hookrightarrow \MA$ associating
with $A\lra \iota_f(C)$ the modal algebra $C$ (and acting similarly on arrows). We call \textit{finite $A$-models} the objects of $A \downarrow \iota_f$.

\begin{theorem}\label{thm:unit}
$\text{PrC}(A) \cong \underleftarrow{\textit{lim}}\, \pi$ and the morphism $A \rightarrow \iota (\underleftarrow{\textit{lim}}\, \pi) \cong \underleftarrow{\textit{lim}} \,\iota \pi$, induced by the universal property of the limit, is the unit of the adjunction ${\rm PrC} \dashv \iota$.
\end{theorem}

\begin{proof}
Let's consider $J$, the full subcategory of $A \downarrow \iota_f$ given by those morphisms $A \rightarrow \iota_f(B)$ that are surjective as functions (such finite $A$-models are called \textit{irreducible}), and the inclusion $J \hookrightarrow (A \downarrow \iota_f)$. Notice that $J$ is a cofiltered category (better, it is a codirected preorder) and that $J$ is initial, so that the limit $\underleftarrow{\textit{lim}} \,\pi$ is isomorphic to the limit of the restriction of $\pi$ via $J$;
from now on, $\pi$ will indicate the restricted cofiltered diagram.

We must show that, for every profinite modal algebra $C$, for every modal algebra $A$
and for every morphism $h: A \lra \iota(C)$ in \MA\ there is unique morphism $\bar h : \underleftarrow{\textit{lim}}\, \pi \rightarrow C$ in ${\rm Pro}\MAf$ such that the composition
$A \xrightarrow{\eta_A} \iota(\underleftarrow{\textit{lim}}\, \pi) \xrightarrow{\iota(\bar h)} \iota(C)$
is $h$ (here $\eta_A$ is the universal map into the limit).
%

Since all objects in ${\rm Pro}\MAf$ are (cofiltered) limits of finite modal algebras, it is easily seen that we can limit the check to the case where $C\in \MAf$. In such a case $h : A \rightarrow \iota(C)$   factorizes trough its image as $A \rightarrow \iota_f(B) \rightarrow \iota_f(C)$, with $A \rightarrow \iota_f(B)$ an object of $J$. Composing the morphism $\eta_A : A \rightarrow \iota (\underleftarrow{\textit{lim}} \,\pi) \cong \underleftarrow{\textit{lim}} \,\iota \pi$, induced by the universal property of the limit, with the image trough $\iota$ of the morphism $\underleftarrow{\textit{lim}}\, \pi \rightarrow \pi(A \rightarrow \iota(B)) = B$ of the limiting cone, we obtain the morphism $A \rightarrow \iota(B)$; composing $\underleftarrow{\textit{lim}} \,\pi \rightarrow B$ with the inclusion $B \rightarrow C$, we can find the desired morphism $\underleftarrow{\textit{lim}}\, \pi \rightarrow C$ in ${\rm Pro}\MAf$ s.t. the  diagram
\begin{center}
\begin{tikzcd}
A \arrow[rd, "\eta_A"'] \arrow[rr, "h"] &                                             & \iota(C) \\
                             & \iota (\underleftarrow{\textit{lim}}\, \pi) \arrow[ru] &
\end{tikzcd}
\end{center}
commutes.
We need to prove that such a morphism is unique. Let's consider two morphisms $f,g : \underleftarrow{\textit{lim}}\, \pi \rightrightarrows C$ in ${\rm Pro}\MAf$ as above; since $C$ is in \MAf, $\underleftarrow{\textit{lim}} \,\pi$ is a cofiltered limit  and the objects of \MAf\
are finitely-copresentable in ${\rm Pro}\MAf$,\footnote{
See Theorem~\ref{teo:charindcompl} and the proof of Theorem~\ref{thm:indcfin} for the dual statement in the dual category \KFrlf.
} $f$ and $g$ factorize, via some $A \rightarrow \iota(B')$ in $J$, as $\underleftarrow{\textit{lim}}\, \pi \rightarrow \pi(A \rightarrow \iota(B')) = B' \overunderset{f'}{g'}{\rightrightarrows} C$ ($J$ is cofiltered, so we can choose the same object for both $f$ and $g$). In the following commutative diagram
\begin{center}
\begin{tikzcd}
A \arrow[rr, "h"] \arrow[d, "\eta_A"'] \arrow[rd]     &                                                                       & \iota(C) \\
\iota (\underleftarrow{\textit{lim}}\, \pi) \arrow[r] & \iota(B') \arrow[ru, "f'", shift left] \arrow[ru, "g'"', shift right] &
\end{tikzcd}
\end{center}
the morphism $A \rightarrow \iota(B')$ is surjective, hence we can conclude that $f' = g'$, implying $f = g$.
\end{proof}

Theorem~\ref{thm:unit} will be useful in the next section to give some concrete description of the left adjoint to the forgetful functor ${\rm Pro}\MAf\lra \SET$.

\section{Transitive Varieties}

All our results easily transfer to finitely approximable varieties of modal algebras (we say that a variety of modal algebras if finitely approximable iff it is generated by its finite members).\footnote{ It does not make sense to consider varieties which are not finitely aprroximable in our context, because profinite algebras are determined by the finite members of a variety.} In fact, the only properties we used about \KFr\ is that \KFrf\  is closed under finite disjoint unions, generated subframes and surjective p-morphisms; since these notions are duals (for finite modal algebras) to the notions of finite products, quotients and subalgebras, it is clear that the required closure properties hold if we consider
 any variety of modal algebras. Hence we have

\begin{theorem}\label{thm:var}
 Let $V$ be a finitely approximable variety of modal algebras and let $V_f$ be the subcategory of finite algebras belonging to $V$. Then the forgetful functor ${\rm Pro}V_f\lra \SET$ is monadic; the unity of the adjointness is given by the limit contruction of Theorem~\ref{thm:unit} applied to free $V$-algebras.
\end{theorem}

Theorem~\ref{thm:var} applies for example to the variety \Int\ of \emph{interior algebras} (also called \emph{S4}-algebras) which is axiomatized by the equations
$$
x\wedge \Diamond x= x \qquad \Diamond \Diamond x =\Diamond x~~.
$$
The class of Kripke frames whose Thomason dual is an interior algebras is formed by \emph{preordered sets}, i.e. by those Kripke frames $(P, \leq)$ where $\leq$ is reflexive and transitive~\cite{misha}.

Another interesting example of finitely approximable variety of modal algebras is given by the Grzegorczyk variety mentioned in the proof of the following result:

\begin{theorem}
 Profinite Heyting algebras are monadic over \SET.
\end{theorem}

\begin{proof}
 Recall that Heyting algebras are distributive lattices with zero and one, endowed with a further binary operation $\to$ such that for all elements $a$ in their support, the posetal functor $a\wedge(-)$ is left adjoint to $a\to(-)$. The category ${\bf Heyt}_f$  of finite Heyting algebras is dual to the category of finite Kripke frames $(P, \leq)$, where $\leq$ is a poset, i.e. it is a reflexive, transitive and antisymmetric relation (morphisms are p-morphisms,  see~\cite{GZ} for a proof).

 Let us now take into consideration the variety  of modal algebras $\bf Grz$; the modal algebras in this variety are the interior algebras satisfying the further axiom
 $$
 \top = \Box(\Box (x\to \Box x)\to x)\to x
 $$
 (recall that $\Box y:= \neg \Diamond \neg y$). The variety is finitely approximable and the class of Kripke frames whose Thomason dual is in $\bf Grz$ is formed by the reflexive and transitive frames $(P, \leq)$ for which there is no infinite ascending chain~\cite{misha}
 $$
 p_0 < p_1 < \cdots < p_i < \cdots
 $$
 (here $p_i< p_{i+1}$ stands for $p_i\leq p_{i+1}$ and $p_{i+1}\neq p_i$).
 Hence a \emph{finite} modal algebra is in $\bf Grz$ iff its dual Kripke frame is a poset. Otherwise said, the categories of finite Heyting algebras and of finite $\bf Grz$-algebras are equivalent. Thus so are the category of profinite Heyting and of profinite $\bf Grz$-algebras. From Theorem~\ref{thm:var}, it follows that the equivalence functor  ${\rm Pro}{\bf Heyt}_f\simeq {\rm Pro}{\bf Grz}_f$ composed with the forgetful functor ${\rm Pro}{\bf Grz}_f\lra \SET$ is monadic (but notice that this is not the forgetful functor
 ${\rm Pro}{\bf Heyt}_f\lra \SET$).
\end{proof}

We conclude the section by an inductive construction of the profinite completion of free interior algebras; the construction (called the `universal model construction'
in terms of Kripke models) exploits Theorem~\ref{thm:unit} and it is well-known in the modal logic literature, where it has been independently introduced for different motivations.\footnote{
In the literature, the construction is supplied only for the case of \emph{finitely generated} free interior algebras, because only in that restricted case the construction enjoys the `definability' properties needed for the applications considered in  the modal logic literature.
}  We follow the exposition of~\cite{ungh}.

Let us call \Pre, \Prelf\  and \Pref\ the categories of preordered sets, of locally finite preordered sets  and of finite preordered sets (morphisms are always p-morphisms); notice that
$(\Prelf)^{op}\simeq {\rm Pro}\Int_f$ and $(\Pref)^{op}\simeq \Int_f$. In locally finite and in finite preordered sets $(P, R)$, there is a notion of \emph{height}    leading to  powerful induction arguments in proofs. The height $h(p_1)$ of a point $p_1$ in a locally finite preordered set $(P, R)$  is defined as the maximum cardinality of chains
$$
 p_1\, R\, p_2\, R \cdots R\, p_n
 $$
such that $p_{i+1} R p_i$ does not hold for $i=1, \dots, n\!-\!1$. Since $(P, R)$ is locally finite and $R$ is transitive, $h(p)$ is a natural number greater or equal to 1.

In the following factorization
\begin{center}
\begin{tikzcd}
{\rm Pro}\Intf \arrow[rr] \arrow[rd, "\iota"', hook] &                        & \SET\\
                                                  & \Int \arrow[ru] &
\end{tikzcd}
\end{center}
the forgetful functor $\Int \rightarrow \SET$ has a left adjoint $F$ (the "free interior algebra" construction); since adjoint functors compose, the left adjoint of ${\rm Pro}\Intf \rightarrow \SET$ associates to a set $X$ the profinite completion of the free interior algebra generated by $X$. We  call this algebra $\text{PrC}_{\bf S4}(F(X))$: according to Theorem~\ref{thm:var}, it  is computed as the limit in ${\rm Pro}\Intf$ of the (filtered) diagram $\pi$ over the category $J$ of irreducible finite $F(X)$-models $F(X) \rightarrow \iota_f(B)$ (the only difference with  respect to the construction of Theorem~\ref{thm:unit} is that now $B$ is a finite interior algebra, not just a finite modal algebra).
From now on, we fix the set $X$ and we just call \emph{finite models} and \emph{finite irreducible  models} the finite and finite  irreducible  $F(X)$-models. It  is easily seen that these are precisely the finite Kripke models (in the standard sense employed in the modal logic literature) over reflexive and transitive Kripke frames,  for the propositional modal language built up from the set of propositional variables $X$.

Via Thomason duality, finite models are (up to iso) of the form $F(X) \rightarrow (\mathcal{P}(P), \Diamond)$, for some $(P,R)\in\Pref$ (here $\Diamond_R$ is the possibility operator associated to the relation $R$, see  the proof sketch of Theorem~\ref{thm:Thomason} above);
according to the duality $(\Pref)^{op}\simeq \Int_f$ morphisms of finite models are commutative triangles
\begin{center}
\begin{tikzcd}
F(X) \arrow[r] \arrow[rd] & {(\mathcal{P}(Q), \Diamond_S)} \arrow[d] \\
                          & {(\mathcal{P}(P), \Diamond_R)}
\end{tikzcd}
\end{center}
induced by the inverse image along some p-morphism $(P,R) \rightarrow (Q,S)$ in $\Pref$. Using that $F$ is left adjoint to the forgetful functor $\Int\ \rightarrow \SET$ and that the power-set functor $\mathcal{P} : \SET \rightarrow \SET^{\text{op}}$ is right adjoint to the power-set functor $\mathcal{P} : \SET^{\text{op}} \rightarrow \SET$, we have the following bijective correspondence
\begin{center}
\begin{tikzcd}[row sep=0.4 em]
                        & F(X) \arrow[r] & {(\mathcal{P}(P), \Diamond_R)} &    & \text{in } \Int  \\
{} \arrow[rrr, no head] &                &                                & {} &                         \\
                        & X \arrow[r]    & \mathcal{P}(P)                 &    & \text{in } \SET \\
{} \arrow[rrr, no head] &                &                                & {} &                         \\
                        & P \arrow[r]    & \mathcal{P}(X)                 &    & \text{in } \SET
\end{tikzcd}
\end{center}
Via this bijection, the category of finite models $(F(X) \downarrow \iota_f)$ is dual to the category having functions $P \rightarrow \mathcal{P}(X)$ (with $(P, R)$ in $\Intf$) as objects and  the p-morphisms $f : (P, R) \rightarrow (Q, S)$ in $\Intf$, s.t. the triangle
\begin{center}
$(\tau)$~~~~~~~~~
\begin{tikzcd}
P \arrow[d] \arrow[r] & \mathcal{P}(X) \\
Q \arrow[ru]          &
\end{tikzcd}
\end{center}
commutes, as morphisms $(P \rightarrow \mathcal{P}(X)) \longrightarrow (Q \rightarrow \mathcal{P}(X))$.

By definition,  a finite model $F(X) \rightarrow \iota_f(B)$ is irreducible iff it is surjective as a function, i.e.  iff  every map
$(F(X) \rightarrow \iota(B'))\longrightarrow (F(X) \rightarrow \iota(B))$ in $(F(X) \downarrow \iota_f)$
induced by  an injective morphism $B'\lra B$ in $\Intf$,
turns out to be an isomorphism. Now noticing that Thomason duality maps injective homomorphisms to surjective p-morphisms, we can conclude that a finite model, seen as a map $P \rightarrow \mathcal{P}(X)$ (for a finite preordered set $(P, R)$) is irreducible iff every p-morphism $f : (P, R) \rightarrow (Q, S)$ commuting the above traiangle $(\tau)$
is an isomorphism. From these observations, the following result follows:

\begin{proposition} Let $(P, R)$ be a finite preordered set. A finite model
  $v : P \rightarrow \mathcal{P}(X)$ is irreducible iff the following two conditions hold for all $p_1, p_2 \in P$,
\begin{itemize}
    \item[a)] $R(p_1) \simeq_{\mathcal{P}(X)} R(p_2) \implies p_1 = p_2$;
    \item[b)] it is not the case that both $R(p_1) \setminus \text{cl}_R(p_1) = R(p_2)$ and $v(\text{cl}_R(p_1)) \subseteq v(\text{cl}_R(p_2))$ hold.
\end{itemize}
Here
 $cl_R(p_1)$ is $\{q\in P\mid p_1Rq~ \&~ qRp_1\}$
and
$R(p_1) \simeq_{\mathcal{P}(X)} R(p_2)$ means~\footnote{
Notice that
we have
$R(p_1)= R^*(p_1)$
in a preordered set.
} that there exists an isomorphism of finite models
\begin{center}
\begin{tikzcd}
R(p_1) \arrow[r, "v|_{R(p_1)}"] \arrow[d, "\cong", "\phi"'] & \mathcal{P}(X) \\
R(p_2) \arrow[ru, "v|_{R(p_2)}"']                  &
\end{tikzcd}
\end{center}
s.t. $\phi(p_1) = p_2$.
\end{proposition}

\begin{proof}
 On one side, if a finite model does not satisfy the above two conditions, it is easy to `reduce' it by taking a suitable coequalizer that identifies $p_1, p_2$
 (use the fact that
 coequalizers of p-morphisms are computed as in \SET,
 see Propositions~\ref{prop:cocomp} and~\ref{prop:cocomplete1}). On the other side, it is not difficult to argue by induction on the cardinality of a finite preordered set $P$ in order to prove that a model
 defined over it and
 satisfying the above two conditions is irreducible.
\end{proof}

From the characterization of irreducible finite models supplied by the above proposition, we can now give an explicit description of the locally finite preordered set dual to the profinite completion of $F(X)$ (this is called the `universal model' over $X$ in the modal logic literature). The construction is by induction, i.e. we successively describe its points of height $n$ by induction on $n$.

For each $n \geq 0$, we define models $v_n : (P_n,R_n) \rightarrow \mathcal{P}(X)$ as follows. $P_0$ is $\emptyset$ and $v_0$ is the empty inclusion. Suppose that $v_n$ has been defined. Put $P_n^+$ equal to the set of pairs $(Y,G)$ such that (i) $Y$ is a nonempty finite subset of $\mathcal{P}(X)$
and $G$ is a generated subframe of $P_n$ comprising at least a point of
height $n$; (ii) if $G = R_n(p)$ for some $p \in P_n$, then $Y \nsubseteq v_n(\text{cl}_{R_n}(p))$.
Put
\begin{align*}
& P_{n+1} :=  P_n \cup \{(y,Y,G)\mid (Y, G)\in P_n^+ ~\&~y\in Y \} \\
& R_{n+1} :=  R_n \cup \{ \langle (y,Y,G), (y', Y, G)\rangle \mid (Y,G)\in P_n^+ ~\&~y, y'\in Y\} \cup
\\
&~~~~~~~~~~~~~~ \cup
\langle \{(y,Y,G), q\rangle \mid (Y,G)\in P_n^+ ~\&~y\in Y~\&~q\in G\} \\
& v_{n+1}(p) :=  \begin{cases}
               v_n(p)       & \quad \text{if } p \in P_n \\
               y            & \quad \text{if } p=(y,Y,G)~\hbox{\rm for some}~(Y,G)\in P^+_n~\&~ y\in Y
               \end{cases}
 \end{align*}
In other words, $(P_n,R_n)$ has been extended by adding all possible ‘not reducible’ bottom clusters of height $n+1$ (a cluster is a subset where the preorder relation is the total relation). It is easily seen by induction that each $v_n : (P_n,R_n) \rightarrow \mathcal{P}(X)$ is irreducible because it satisfies a) and b). Moreover, each irreducible model can be $\mathcal{P}(X)$-embedded into $v_n : (P_n,R_n) \rightarrow \mathcal{P}(X)$ for some $n$. Hence the chain formed by these models is final in the colimit construction that gives
the Kripke frame dual to $\text{\rm PrC}_{\bf S4}(F(X))$,
so
such a Kripke frame is the set-theoretic union on this chain.
To conclude, we have that
\begin{theorem}
$\text{\rm PrC}_{\bf S4}(F(X)) ~~\cong~~ \langle \mathcal{P}( \bigcup_{n \geq 0} P_n), \Diamond_{\bigcup R_n}\rangle~~~.$
\end{theorem}

\section{Conclusions}

We exploited Thomason duality in order to identify a full subcategory  of the category of Kripke frames  which is dual to profinite modal algebras; this subcategory (the category of locally finite Kripke frames \KFrlf)  is surprisingly well-behaved and the forgetful functor
from it
into the category of sets turns out to be
comonadic.
This  makes profinite modal algebras monadic over \SET\
and extends to modal logic the well-known monadicity result for complete atomic Boolean algebras. The monad itself is interesting as,
in terms of dual Kripke models,
it corresponds to the universal model construction well-known to the modal logic community.

Further work is needed to investigate the exactness property of categories like \KFrlf, especially in correspondence to  transitive subvarieties of \MA. From preliminary analysis, it seems that regularity can  be obtained in some limited but relevant cases, whereas exactness could
incur in couterexamples similar to those
exhibited in~\cite{GZ} for the opposite of the category of finitely presented Heyting algebras.

Finally, from the results of Section~\ref{sec:kfrlf}, it follows that \KFrlf\ is a locally finitely presentable category, hence it is the category of models of an essentially algebraic first-order theory: it would be interesting to identify such a theory (the task is made difficult by the fact that morphisms in \KFrlf\ are p-morphisms, not just stable maps).

\bibliographystyle{unsrt}
\bibliography{biblio}  






\section{Appendix}

We supply here some missing proofs of background results. Recall that an \emph{atom} in a Boolean algebra $B$ is a non-zero minimal element; if the algebra is complete, atoms can be equivalently defined as Join-irreducible elements, i.e. as the elements $a$ such that, for any $C\subseteq B$, we have that $a\leq \bigvee C$ implies that there is
$c\in C$ such that $a\leq c$.
\vskip 2mm
\noindent
\textbf{Proposition~\ref{prop:caba}}
\emph{
 The following conditions are equivalent for a Boolean algebra $B$:
  \begin{itemize}
   \item $B$ is complete and atomic;
   \item $B$ is isomorphic to a powerset Boolean algebra;
   \item $B$ is complete and satisfies the infinitary distributive law:
     \begin{align*}
      \bigwedge_{i \in I} \bigvee X_i = \bigvee_{f \in \prod_{i \in I} X_i} \bigwedge_{i \in I} f(i)
     \end{align*}
     for every family $\{X_i\}_{i \in I}$ of subsets of $B$.
  \end{itemize}
  }
\begin{proof}
A complete Boolean algebra $B$ is atomic iff each of its elements $b \in B$ can be written as the disjunction of the atoms $a\in B$ such that $a\leq b$.
For a Boolean algebra $B$, we can consider the  morphism
\begin{equation}\label{eq:eta}
\eta_B : B \rightarrow \mathcal{P}(\Atm(B)),
\end{equation}
sending each element $b \in B$ to the set of atoms below it; by the observation above, if $B$ is complete and atomic, $\eta_B$ must be an isomorphism.

If $B$ is isomorphic to a powerset Boolean algebra, then it's easy to see that it is complete and satisfies the infinitary distributive law.

If $B$ is complete and satisfies the infinitary distributive law, then we can write
     \begin{align*}
      \top = \bigwedge_{b \in B} (b \vee \neg b) = \bigvee_{f \in \prod_{b \in B} \{b, \neg b\}} \bigwedge_{b \in B} f(b)
     \end{align*}
We define $a(f) := \bigwedge_{b \in B} f(b)$ for all $f \in \prod_{b \in B} \{b, \neg b\}$. Given $b' \in B$ and $f \in \prod_{b \in B} \{b, \neg b\}$ s.t. $a(f) \neq \bot$,
    \begin{align*}
    b' \wedge a(f) \neq \bot \implies f(b') = b' \text{ (by definition of } a(f) \text{)} \implies a(f) \leq b'
    \end{align*}
and the viceversa is obviously true; this implies that $\{a(f) \; | \; f \in \prod_{b \in B} \{b, \neg b\} ~\& \\ a(f) \neq \bot\}$ is the set of atoms of $B$ and, since
     \begin{align*}
      b' = b' \wedge \top = b' \wedge \bigvee_{f \in \prod_{b \in B} \{b, \neg b\}} a(f) = \bigvee_{f \in \prod_{b \in B} \{b, \neg b\}} b' \wedge a(f)
     \end{align*}
for every $b' \in B$, we conclude that $B$ is atomic.
\end{proof}

\noindent
\textbf{Theorem~\ref{thm:Tarski}}
\emph{
[Tarski duality] \CABA\ is dual to \SET.
}
\begin{proof}
We have a functor
    \begin{center}
    \begin{tikzcd}
    \mathcal{P} &[-3em] : &[-3em] \SET^{\text{op}} \arrow[r] & \CABA
    \end{tikzcd}
    \end{center}
sending a set $X$ into the Boolean algebra $\mathcal{P}(X)$, which we proved to be complete and atomic. Given a function $f : X \rightarrow X'$, $\mathcal{P}(f) : \mathcal{P}(X') \rightarrow \mathcal{P}(X)$ sends a subset $S' \subseteq X'$ into the pre-image $f^{-1}(S')$; it's easy to see that it defines a morphism of complete atomic Boolean algebras.

Viceversa, we have the functor
\begin{center}
\begin{tikzcd}
\Atm &[-3em] : &[-3em] \CABA \arrow[r] & \SET^{\text{op}}
\end{tikzcd}
\end{center}
sending a complete atomic Boolean algebra $B$ into the set of its atoms $\Atm(B)$. Given a morphism $\mu : B \rightarrow B'$ of complete atomic Boolean algebras, we can take its left adjoint $\mu^* : B' \rightarrow B$. $\mu^*$ restricts to a function $\Atm(\mu) : \Atm(B') \rightarrow \Atm(B)$: in fact, given $a' \in \Atm(B')$ and $C \subseteq B$, if $\mu^*(a') \leq \bigvee C$ (i.e. $a' \leq \mu \left( \bigvee C \right) = \bigvee \mu(C)$), then there exists $c \in C$ s.t. $a' \leq \mu(c)$, i.e. there exists $c \in C$ s.t. $\mu^*(a') \leq c$.

The isomorphism~\eqref{eq:eta} of complete atomic Boolean algebras $\eta_B : B \rightarrow \mathcal{P}(\Atm(B))$  is natural in $B$ and the function $\epsilon_X : X \rightarrow \Atm(\mathcal{P}(X))$, sending $x \in X$ to the singleton $\{x\}$, is an isomorphism in \SET\ and it's natural in $X$. These two maps are the unity and counity of the equivalence.
\end{proof}

\noindent
\textbf{Theorem~\ref{thm:Thomason}}
\emph{
 [Thomason duality] \CAMA\ is dual to \KFr.
 }
\begin{proof}
Similarly to Tarski duality, we have a pair of functors
\begin{center}
\begin{tikzcd}
\mathcal{P} &[-3em] : &[-3em] \KFr^{\text{op}} \arrow[r] & \CAMA
\end{tikzcd}
\end{center}
and
\begin{center}
\begin{tikzcd}
\Atm &[-3em] : &[-3em] \CAMA \arrow[r] & \KFr^{\text{op}}
\end{tikzcd}
\end{center}
In fact, given a Kripke frame $(P,R)$, we can define a diamond operator $\Diamond_R : \mathcal{P}(P) \rightarrow \mathcal{P}(P)$ sending $A \in \mathcal{P}(P)$ to $\{p \in P \; | \; \exists p' \in A \text{ s.t. } pRp' \}$ (it's easy to check that this operator commutes with arbitrary joins) making $\mathcal{P}(f) : \mathcal{P}(Q) \rightarrow \mathcal{P}(P)$ a morphism of completely-additive complete atomic modal algebras for every p-morphism of Kripke frames $f : P \rightarrow Q$.

Viceversa, given a completely-additive complete atomic modal algebra $M$, we can define the following relation on the set of its atoms: for $a,a' \in \Atm(M)$,
     \begin{align*}
      a R_M a' \iff a \leq \Diamond a'
     \end{align*}
For every morphism $\mu : M \rightarrow N$ in \CAMA, the function $\Atm(\mu) : \Atm(N) \rightarrow \Atm(M)$ is stable and open wrt the relations $R_N$ and $R_M$:
\begin{itemize}
    \item[$\bullet$] if $b R_N b'$ (i.e. $b \leq \Diamond b'$), then $\mu^*(b) \leq \mu^*(\Diamond b') \leq \Diamond \mu^*(b')$ (from $b' \leq \mu(\mu^*(b'))$, we get $\Diamond b' \leq \Diamond \mu(\mu^*(b')) = \mu(\Diamond \mu^*(b'))$, i.e. $\mu^*(\Diamond b') \leq \Diamond \mu^*(b')$); this means  $\mu^*(b) R_M \mu^*(b')$;
    \item[$\bullet$] if $\mu^*(b) R_M a'$, i.e. $\mu^*(b) \leq \Diamond a'$, then  $b \leq \mu(\Diamond a') = \Diamond \mu(a') = \Diamond \bigvee \{b' \; | \; b' \in \Atm(N) \text{ and } b' \leq \mu(a')\} = \bigvee \{\Diamond b' \; | \; b' \in \Atm(N) \text{ and } b' \leq \mu(a')\}$ (being $N$ atomic and completely-additive); but then there exists $b' \in \Atm(N)$ s.t. $b \leq \Diamond b'$ and $b' \leq \mu(a')$ (i.e. $\mu^*(b') \leq a'$, hence $\mu^*(b') = a'$, being $a'$ and $\mu^*(b')$  atoms), that is there exists $b' \in \Atm(N)$ s.t. $b R_N b'$ and $\mu^*(b') = a'$.
\end{itemize}

As for Tarski duality, we have a pair of natural isomorphisms, $\eta_M : M \rightarrow \mathcal{P}(\Atm(M))$ in \CAMA\ and $\epsilon_P : P \rightarrow \Atm(\mathcal{P}(P))$ in \KFr\ satisfying triangular identities.
\end{proof}

\end{document}